\magnification=\magstep1

\hsize=15 true cm  
\vsize=23 true cm
\baselineskip13pt 

\font\footfont=cmr8
\font\bigbf=cmbx12
\font\csc=cmcsc10 
\font\smallrm=cmr8 
\font\cyr=wncyr10 
\def\j3{{\rm\u{\cyr i}}}

\def\square{{\ \vrule height0.5em width0.5em depth-0.0em}}
\def\D{\mathop{\rm Var}}
\def\MSE{{\rm MSE}}
\def\MISE{{\rm MISE}}

\def\MSE{{\rm MSE}}

\def\mv{{\rm v.p.}}
\def\half{\hbox{$1\over2$}} 
\def\V{{\bf V}}
\def\E{{\rm E}}
\def\d{{\rm d}}
\def\section{\bigskip}
\def\start{\medskip\noindent} 
 
\def\RR{\mathord{I\kern-.3em R}}

\centerline{\bigbf Upper Bounds for the I-MSE and max-MSE} 
\centerline{\bigbf of Kernel Density Estimators} 

\bigskip
\centerline{\bf 
Nils Lid Hjort$^1$ and  
Nikolai G.~Ushakov$^{2,}$\footnote{{\rm *}}{\footfont Partially 
supported by RFBR Grants 97-01-00273,
98-01-00621 and 98-01-00926, and by INTAS-RFBR Grant IR-97-0537}}

\medskip
\centerline{\bf December 1999} 
\medskip 
\noindent $^1$Department of Mathematics, University of Oslo, Norway

\noindent $^2$Russian Academy of Sciences, 
Chernogolovka, Russia

\medskip
{\smallskip\narrower\noindent\baselineskip12pt 
{\csc Abstract.} 
The performance of kernel density estimators is usually 
studied via Taylor expansions and asymptotic approximation arguments,
in which the bandwidth parameter tends to zero with increasing sample size.
In contrast, this paper focusses directly on the finite-sample situation.
Informative upper bounds are derived both for 
the integrated and the maximal mean squared error function. 
Results are reached for the traditional case,  
where the kernel is a probability density function, 
under various sets of assumptions on the underlying density
to be estimated. Results are also derived for the 
important non-conventional case of the sinc kernel,
which is not integrable and also takes negative values. 
We pin-point ways in which the sinc-based estimator 
performs better than the conventional kernel estimators. 
When proving our results we rely on methods related 
to characteristic and empirical characteristic functions. 

\smallskip\noindent
{\csc Key words:} \sl
characteristic functions, 
density estimation,  
finite-sample performance, 
max-MSE, 
sinc kernel, 
upper bounds
\smallskip}

\section 
\centerline{\bf 1. Introduction}

\start 
In this article we derive some rigorous upper bounds 
for the estimation error of kernel density estimators for finite 
values of the sample size $n$, in terms of choices of 
the kernel function $K$ and the bandwidth $h=h_n$. 
These bounds are by construction non-asymptotic,
and are useful when one needs to secure a certain precision 
of an estimate for a given (finite) value of $n$,
for broad classes of densities. 
We study both smooth cases (where the density to be estimated 
is one or more times differentiable) 
and non-smooth cases (the underlying density function 
is not supposed to be differentiable or even continuous).
The machinery of characteristic and empirical characteristic
functions is used, and relevant general results are 
established in Section 2. 

In Section 3 conventional kernel estimators will be considered,
i.e.~estimators whose kernels are probability density
functions. These estimators always produce estimates which
are densities. 
We term a kernel density estimator non-conventional if its 
kernel function is not a probability density, i.e.~it may 
take negative values or/and does not integrate to one 
(or even is not integrable). 
Such non-conventional estimators are studied in Section 4,
with particular attention to the sinc kernel;
see also Glad, Hjort and Ushakov (1999a). 
Such estimators, based on higher order kernels,
superkernels or the sinc kernel, often provide better estimation
precision, but have an essential disadvantage: they produce 
estimates which are not probability density functions, 
i.e.~may take negative values or/and do not integrate to one. 
However, this defect can be corrected afterwards 
without loss of their performance properties
(see Glad, Hjort and Ushakov, 1999b). 

A discussion of our results, with a view towards 
their use in density estimation problems, 
is given in the final Section 5. Topics there 
include new strategies for bandwidth selection. 

\section 
\centerline{\bf 2. Auxiliary results, via characteristic functions} 

\start
In this paper we use the characteristic function approach
to studying performance of density estimators, rather than
the traditional Taylor expansions and asymptotic approximations. 
Therefore we first express some basic concepts of kernel density estimators 
in terms of characteristic functions. 

Let $X_1,\ldots,X_n$ be independent and
identically distributed random variables with absolutely continuous
distribution function $F(x)$, density function $p(x)$, and
characteristic function $f(t)$. The kernel density estimator
associated with the sample $X_1,\ldots,X_n$ is defined as
$$p_n(x)=p_{n,h}(x)=
{1\over{n}}\sum_{j=1}^n K_h(x-X_j),\eqno(2.1)$$
where $K(x)$ is the kernel function with scaled version
$K_h(x)=h^{-1}K(h^{-1}x)$ and $h=h_n$ is a positive number 
(depending on $n$) called the bandwidth or the smoothing parameter.
We do not necessarily demand that $K$ is integrable 
(sometimes the best estimators correspond to nonintegrable kernels). 
However, we suppose that $K$ is square integrable,
and in addition that it is integrable in the sense
of the Cauchy principal value with 
$\mv\int_{-\infty}^\infty K(x)\,\d x=1$, 
in which 
$$\mv\int_{-\infty}^\infty=\lim_{T\to\infty}
  \lim_{\epsilon\to0}\Bigl[\int_{-T}^{-\epsilon}+\int_{\epsilon}^{T}\Bigr].$$
Under these assumptions the Fourier transform of $K$ can be defined as
$$\varphi(t)=\mv\int_{-\infty}^\infty e^{itx}K(x)\,\d x$$
(see Chapter 4 of Titchmarsh, 1937). 
In the following we will omit integration limits when
the integral is to be taken over the full real line. 

Let $\hat p_n$ be an estimator (not necessarily a kernel estimator) 
of $p$ associated with the sample $X_1,\ldots,X_n$. The bias,
the mean squared error (MSE) and the mean integrated squared error (MISE) of
$\hat p_n$ are defined, respectively, as
$$B_n(\hat p_n(x))=\E \hat p_n(x)-p(x), $$
$$\MSE(\hat p_n(x))=\E\{\hat p_n(x)-p(x)\}^2, $$
and
$$\MISE(\hat p_n)=\int\MSE(\hat p_n(x))\,\d x
  =\E \int\{\hat p_n(x)-p(x)\}^2\,\d x. \eqno(2.2)$$
In case of the kernel estimator $p_n$, defined by (2.1), 
the bias may be expressed via the convolution as 
$$B_n(p_n(x))=(K_h\star p)(x)-p(x)=\int K_h(x-y)p(y)\,\d y-p(x). $$
Since convolution is a kind of smoothing, the bias of the kernel
estimator is the difference between a smoothed density and the density itself.
The mean squared error 
admits a well-known decomposition into variance and squared bias,
with consequent $\MISE$ representation 
$$\MISE(\hat p_n)=\int B_n^2(\hat p_n(x))\,\d x+
  \int \D(\hat p_n(x))\,\d x.$$

Note that together with $\MSE$ and $\MISE$ other measures of 
deviation may be used. Among them, the mean absolute error
$\E |\hat p_n(x)-p(x)|$ and its integral 
are especially important (see Devroye and Gy\"orfi, 1985). 
In the present article we restrict attention to $\MSE$ and $\MISE$,
however. 

For a real valued function $g$ we will use the following notation, 
provided the integrals exist: 
$$\mu_k(g)=\int |x|^kg(x)\,\d x
  \quad {\rm and} \quad R(g)=\int g^2(x)\,\d x. $$
If the kernel $K$ is a probability density function, 
and the density to be estimated 
is twice differentiable and with square integrable second order derivative,
then it is well known that the best order of estimation accuracy 
in terms of MISE is $O(n^{-4/5})$; see also Section 5.1. 
However, if we permit the kernel not to be a density, 
then the order can be improved. For example,  
if $p$ is the normal density and $K$ is the sinc kernel, 
i.e.~$K(x)=\sin x/(\pi x)$, then
$$\min_{h>0}\MISE(p_n)=O\Bigl({{\sqrt{\log n}}\over n}\Bigr)
  \quad {\rm as\ }n\to\infty\,; $$
see Section 4 below and Glad, Hjort and Ushakov (1999a). 

We now express basic characteristics of density estimators in terms
of Fourier transforms and establish some auxiliary results.

Let $\hat f_n$ denote the Fourier transform of an estimator 
$\hat p_n$. Making use of the inversion formula for densities
and the Parseval--Plancherel identity 
we easily obtain the following formulae: 
$$B_n(\hat p_n(x))={1\over{2\pi}}\int e^{-itx}
  \{\E \hat f_n(t)-f(t)\}\,\d t,\eqno(2.3)$$
$$\eqalign{
\MSE(\hat p_n(x))
&=\E\Bigl\{{1\over{2\pi}}\int e^{-itx}[\hat f_n(t)-f(t)]\,\d t
  \Bigr\}^2 \cr 
&={1\over{(2\pi)^2}}\int\int e^{-i(u+v)x}
  \E \{(\hat f_n(u)-f(u))(\hat f_n(v)-f(v))\}\,\d u\,\d v, \cr}\eqno(2.4)$$
and 
$$\MISE(\hat p_n)={1\over{2\pi}}
  \int \E |\hat f_n(t)-f(t)|^2\,\d t.\eqno(2.5)$$

In the remainder of this section, we will consider only kernel
estimators and suppose that the kernel $K$ is a probability
density function, i.e.~it is nonnegative and integrates to one.
Study the empirical characteristic function associated with sample
$X_1,\ldots,X_n$, 
$$f_n(t)={1\over n}\sum_{j=1}^ne^{itX_j}.$$
The characteristic function of the estimator $p_{n,h}(x)$ is
$f_n(t)\varphi(ht)$, where $\varphi(t)=\int e^{itu}K(u)\,\d u$
is the characteristic function of the kernel. 
And the kernel estimator (2.1) can be expressed in terms of $f_n$ as 
$$p_{n,h}(x)={1\over{2\pi}}
  \int e^{-itx}f_n(t)\varphi(ht)\,\d t.$$
Now, taking into account that 
$$\E f_n(u)f_n(v)=(1-1/n)f(u)f(v)+(1/n)f(u+v) $$
and
$$\E |f_n(t)|^2=(1-1/n)|f(t)|^2+1/n,$$
we can write (2.3)--(2.5) in the form
$$B_n(p_n(x))={1\over{2\pi}}
  \int e^{-itx}f(t)\{\varphi(ht)-1\}\,\d t,\eqno(2.6)$$
$$\eqalign{
\MSE(p_n(x))=
&{1\over{(2\pi)^2}}
  \int\int e^{-i(u+v)x}\Big[{1\over n}\varphi(hu)\varphi(hv)f(u+v) \cr 
&+\Big\{\Big(1-{1\over n}\Big)\varphi(hu)\varphi(hv)
  -2\varphi(hu)+1\Big\}f(u)f(v)
  \Big]\,\d u\,\d v \cr}\eqno(2.7)$$
and
$$\MISE(p_n)={1\over{2\pi}}\Big[
  \int |f(t)|^2|1-\varphi(ht)|^2\,\d t
  +{1\over n}\int|\varphi(ht)|^2\{1-|f(t)|^2\}\,\d t\Big].\eqno(2.8)$$

From (2.6) we immediately obtain
$$|B_n(p_n(x))|\le{1\over{2\pi}}\int |f(t)|
  \,|1-\varphi(ht)|\,\d t.\eqno(2.9)$$

\smallskip
{\bf Lemma 1.} \sl For each $x$, 
$$\MSE(p_n(x))
  \le\Bigl\{{1\over{2\pi}}\int |f(t)|
  \,|1-\varphi(ht)|\,\d t\Bigr\}^2
  +{a(x)\over{\pi nh}}\int|\varphi(t)|\,\d t, \eqno(2.10)$$
where $a(x)=(K_h\star p)(x)$. If $p$ is bounded by $a$, then
$$\sup_x\MSE(p_n(x))
  \le\Bigl\{{1\over{2\pi}}\int |f(t)|
  \,|1-\varphi(ht)|\,\d t\Bigr\}^2
  +{a\over{\pi nh}}\int|\varphi(t)|\,\d t. $$
\rm

\smallskip
{\bf Proof.} It suffices to prove the first statement,
since $a(x)=\int p(x-y)K_h(y)\,\d y\allowbreak \le a$ for all $x$. 
Making use of relation (2.7), we obtain
$$\eqalign{\MSE(p_n(x))
=&\Bigl[{1\over{2\pi}}\int e^{-itx}f(t)
  \{1-\varphi(ht)\}\,\d t\Bigr]^2 \cr 
&+{1\over n}{1\over{(2\pi)^2}}\int \int 
  e^{-i(u+v)x}\varphi(hu)\varphi(hv)f(u+v)\,\d u\,\d v \cr 
&-{1\over n}{1\over{(2\pi)^2}}\int \int 
  e^{-i(u+v)x}\varphi(hu)\varphi(hv)f(u)f(v)\,\d u\,\d v. \cr}$$
The first term on the right hand side is dominated by the first term
of the right hand side of (2.10). Let us then estimate the absolute value of
the second (denoted by $T_2$) and third (denoted by $T_3$) terms.
We have
$$T_2={1\over n}{1\over{2\pi}}\int \varphi(hu)
  \Bigl\{{1\over{2\pi}}\int 
  e^{-i(u+v)x}\varphi(hv)f(u+v)\,\d v\Bigr\}\,\d u.$$
The term in brackets, being transformed to the form
$${1\over{2\pi}}\int 
e^{-itx}\varphi(h(t-u))f(t)\,\d t,$$
is equal to $\int p(x-y)K_h(y)e^{-iuy}\,\d y$ 
(since $\varphi(h(t-u))f(t)$ is the Fourier transform of the convolution
of functions $p(x)$ and $K_h(x)e^{-iux}$),
and clearly 
$$\Bigl|\int  p(x-y)K_h(y)e^{-iuy}\,\d y\Bigr|\le a(x). $$
Hence
$$|T_2|\le{a(x)\over n}{1\over{2\pi}}\int|\varphi(ht)|\,\d t.$$
Furthermore,
$$\eqalign{|T_3|=
&{1\over n}\Bigl|
  {1\over{2\pi}}\int  e^{-iux}\varphi(hu)f(u)\,\d u\,
  {1\over{2\pi}}\int  e^{-ivx}f(v)\varphi(hv)\,\d v\Bigr| \cr 
&\le {1\over n}(K_h\star p)(x){1\over 2\pi}\int |f(v)|\,|\varphi(hv)|\,\d v 
 \le {a(x)\over n\,2\pi}\int|\varphi(hv)|\,\d v. \cr}$$
Thus we finally obtain (2.10). $\square$

\smallskip
{\bf Lemma 2.} \sl
$$\MISE(p_{n,h})\le{1\over{2\pi}}
  \Bigl\{\int |f(t)|^2\,|1-\varphi(ht)|^2\,\d t+
  {1\over {nh}}\int |\varphi(t)|^2\,\d t\Bigr\}. $$
\rm

\smallskip
This lemma immediately follows from relation (2.8).

We conclude this section with some inequalities for 
characteristic functions and which will be used below.

\smallskip
{\bf Lemma 3.} \sl Let $F$ be a distribution function with
characteristic function $f$. If the first order absolute
moment $\beta_1=\int |x|\,\d F(x)$ is finite, then
$$|1-f(t)|\le\beta_1|t| \quad {\sl for\ all\ real\ }t. $$
If $F$ has null expectation and finite variance $\sigma^2$, then
$$|1-f(t)|\le \half\sigma^2t^2\quad {\sl for\ all\ real\ }t. $$
\rm

\smallskip
{\bf Proof.} Observe that for any positive integer $n$ and any $x>0$,
$$\Bigl|e^{ix}-1-{{ix}\over{1!}}-\ldots-{{(ix)^{n-1}}\over{(n-1)!}}
\Bigr|\le{{x^n}\over{n!}}\eqno(2.11)$$
(see for example Feller, 1971, Chapter 15).
The first inequality of the lemma follows quickly via 
$$|1-f(t)|\le\int |1-e^{itx}|\,\d F(x)
  \le\int |tx|\,\d F(x)=\beta_1|t|.$$
To prove the second inequality, 
we obtain, again making use of (2.11), 
$$\eqalign{|1-f(t)|=
&\Bigl|\int \{e^{itx}-1\}\,\d F(x)
 \Bigr|=\Bigl|\int (e^{itx}-1-itx)\,\d F(x)\Bigr| \cr 
&\le\int |e^{itx}-1-itx|\,\d F(x)
 \le\int \half t^2x^2\,\d F(x)=\half\sigma^2t^2. \cr}$$
Along the same lines one may prove for example that 
$|f(t)-(1-\half\sigma^2t^2)|\le {1\over6}|t|^3\int|x|^3\,\d F(x)$. 
$\square$

\smallskip
Let $g$ be a real-valued function defined on an interval
$[a,b]$ of the real line.
The total variation of $g$ on $[a,b]$ is defined as
$$\V_a^b(g)=\sup\sum_{i=1}^n|g(x_i)-g(x_{i-1})|$$
where the supremum is taken over all $n$ and all collections 
$x_0,\ldots,x_n$ such that $a=x_0<\cdots<x_n=b$. The total
variation on the whole real line is defined as
$$\V_{-\infty}^\infty(g)=\lim_{x\to\infty}\V_{-x}^x(g).$$
In the case $\V_{-\infty}^\infty(g)$ we omit limits and
write $\V(g)$.
A function $g$ is said to be a function of bounded total
variation if $\V(g)<\infty$ 
(or $\V_a^b(g)<\infty$ if it is considered on an interval $[a,b]$). 
Note that if $g$ has an integrable derivative, then 
$\V_a^b(g)=\int_a^b |g'|\,\d x$. 

\smallskip
{\bf Lemma 4.} \sl Let $p$ be a probability density
and $f$ the corresponding characteristic function. If
$p$ is $m-1$ times differentiable, 
and $p^{(m-1)}$ is a function of bounded variation, then
$$|f(t)|\le \V(p^{(m-1)})/|t|^m $$
for all real $t$ (by definition, $p^{(0)}=p$).
\rm

\smallskip
A proof of this lemma is contained in Ushakov and Ushakov (1999).

\section 
\centerline{\bf 3. Density estimators with conventional kernels}

\start 
First we study the `smooth' case, i.e.~when the density
to be estimated is one or several times differentiable.

\smallskip
{\bf Theorem 1.} \sl Let $p$ be twice differentiable, 
with $p''$ a function of bounded variation, $\V(p'')=V_2<\infty$. 
If the kernel $K$ has null expectation,
and $h_n=h_0n^{-1/5}$ $(h_0$ being some constant$)$, then
$$\MISE(p_n)\le
  \Bigl\{
  {{3\mu_2^2(K)V_2^{5/3}}\over{10\pi}}h_0^4+{{R(K)}\over{h_0}}\Bigr\}
  n^{-4/5}.\eqno(3.1)$$
\rm

\smallskip
{\bf Proof.} Due to Lemma 4, we have
$$|f(t)|\le\cases{
1&for $|t|\le V_2^{1/3}$,\cr
V_2/|t|^3&for $|t|>V_2^{1/3}$,\cr}\eqno(3.2)$$
and, due to Lemma 3,
$|1-\varphi(h_nt)|\le \half \mu_2(K)h_n^2t^2$ for all $t$. Hence
$$\eqalign{
\int |f(t)|^2|1-\varphi(h_nt)|^2\,\d t
&\le \half \mu_2^2(K)h_n^4\int_0^{V_2^{1/3}}t^4\,\d t
  +\half\mu_2^2(K)h_n^4V_2^2\int_{V_2^{1/3}}^\infty
  {{\,\d t}\over{t^2}} \cr
&=(3/5)\mu_2^2(K)V_2^{5/3}h_0^4\,n^{-4/5}. \cr}\eqno(3.3)$$
Further, using the Parseval--Plancherel identity, we get
$${1\over {nh_n}}{1\over{2\pi}}\int|\varphi(t)|^2\,\d t
  ={1\over{h_0}}n^{-4/5}\int  K^2(x)\,\d x
  ={{R(K)}\over{h_0}}n^{-4/5}.\eqno(3.4)$$
From (3.3), (3.4) and Lemma 2, we obtain (3.1).
$\square$ 


\smallskip
{\bf Corollary.} \sl Let the conditions of Theorem 1 be satisfied.
Then for each $n$, 
$$\min_{h>0}\MISE(p_n)\le
  \Bigl({{3\cdot 5^4}\over{2^9\pi}}\Bigr)^{1/5}
  \{\mu_2^2(K)R^4(K)\}^{1/5}V_2^{1/3}n^{-4/5}, $$
with minimum of the upper bound attained for 
$$h_n=\Bigl\{{5\pi\over 6}{R(K)\over \mu_2^2(K)}\Bigr\}^{1/5}
  V_2^{-1/3}n^{-1/5}. $$
\rm

\smallskip\noindent 

If $p$ is only one time differentiable or/and the expectation of
$K$ does not equal zero, then results are weaker.

\smallskip
{\bf Theorem 2.} \sl Let $p$ be differentiable with 
$p'$ a function of bounded variation, 
$\V(p')=V_1<\infty$. 
If $h_n=h_0n^{-1/3}$, then
$$\MISE(p_n)\le
  \Bigl\{{4\over{3\pi}}\mu_1^2(K)V_1^{3/2}h_0^2+{{R(K)}\over{h_0}}
  \Bigr\}n^{-2/3}.\eqno(3.5)$$
\rm

\smallskip
{\bf Proof.} Due to Lemmas 3 and 4,
$$|f(t)|\le\cases{
1&for $|t|\le V_1^{1/2}$,\cr
V_1/|t|^2&for $|t|>V_1^{1/2}$,\cr}$$
and $|1-\varphi(h_nt)|\le\mu_1(K)h_n|t|$ for all $t$. 
Hence (see the proof of Theorem 1),
$$\int |f(t)|^2|1-\varphi(h_nt)|^2\,\d t
\le (8/3)\mu_1^2(K)V_1^{3/2}h_0^2n^{-2/3}.\eqno(3.6)$$
And, as in the proof of Theorem 1,
$${1\over{2\pi}}{1\over{nh_n}}\int 
|\varphi(t)|^2\,\d t={{R(K)}\over{h_0}}n^{-2/3}.\eqno(3.7)$$
From (3.6), (3.7) and Lemma 2, we obtain (3.5).
$\square$
\eject 

\smallskip
{\bf Corollary.} \sl Let the conditions of Theorem 2 be satisfied. Then
for each $n$,
$$\min_{h>0}\MISE(p_n)
  \le (9/\pi)^{1/3}\mu_1^{2/3}(K)\sqrt{V_1}R(K)^{2/3}\,n^{-2/3}. $$
\rm

\smallskip
Theorems 1 and 2 give bounds for the integral deviation of the 
mean squared error of a kernel estimator from zero. Now we obtain
bounds for the sup deviation, in terms of 
$$A(K)={1\over{2\pi}}\int |\varphi(t)|\,\d t.$$

\smallskip
{\bf Theorem 3.} \sl Let $p$ be three times differentiable
with $p'''$ a function of bounded variation, 
$\V(p''')=V_3<\infty$, and let $p$ be bounded by $a$. 
If $h_n=h_0n^{-1/5}$, then
$$\sup_x\MSE(p_n(x))\le
  \Bigl\{{{4}\over{9\pi^2}}\mu_2^2(K)V_3^{3/2}h_0^4
  +{{2aA(K)}\over{h_0}}\Bigr\}\,n^{-4/5}.$$
\rm

\smallskip
{\bf Proof.} Due to Lemma 4,
$$|f(t)|\le\cases{
1&for $|t|\le V_3^{1/4}$,\cr
V_3/|t|^4&for $|t|>V_3^{1/4}$,\cr}$$
and, due to Lemma 3,
$|1-\varphi(h_nt)|\le \half\mu_2(K)h_n^2t^2$ for all $t$. Hence
$$\eqalign{
\int |f(t)|\,|1-\varphi(h_nt)|\,\d t
&\le\mu_2(K)h_n^2\Bigl(\int_0^{V_3^{1/4}}t^2\,\d t
  +V_3\int_{V_3^{1/4}}^\infty {{\,\d t}\over{t^2}}\Bigr) \cr 
&=(4/3)\mu_2(K)V_3^{3/4}h_0^2n^{-2/5}. \cr}$$
To get the result it suffices now to apply Lemma 1.
$\square$

\smallskip
{\bf Corollary.} \sl Let the conditions of Theorem 3 be satisfied. 
Then for each $n$,
$$\min_{h>0}\sup_x\MSE(p_n(x))
  \le 5(36\pi^2)^{-1/5}\mu_2^{2/5}(K)V_3^{3/10}A^{4/5}(K)a^{4/5}\,n^{-4/5}, $$
with minimum of the upper bound being attained for 
$$h_n=\Bigl\{{9\pi^2\over 8}{A(K)\over \mu_2^2(K)}\Bigr\}^{1/5}
  a^{1/5}V_3^{-3/10}\,n^{-1/5}. $$
\rm

\smallskip
{\bf Theorem 4.} \sl Let $p$ be twice differentiable
with $p''$ a function of bounded variation, and 
let $p$ be bounded by $a$. If $h_n=h_0n^{-1/3}$, then
$$\sup_x\MSE(p_n(x))\le
\Bigl\{
{{9}\over{4\pi^2}}\mu_1^2(K)V_2^{4/3}h_0^2+{{2aA(K)}\over{h_0}}\Bigr\}
n^{-2/3}.\eqno(3.8)$$
\rm

\smallskip
{\bf Proof.} Using (3.2) and the second inequality of Lemma 3, we have 
$|1-\varphi(h_nt)|\le\mu_1(K)h_n|t|$. This leads to 
$$\eqalign{
\int |f(t)|\,|1-\varphi(h_nt)|\,\d t
&\le 2\mu_1(K)h_n\int_0^{V_2^{1/3}}t\,\d t+2V_2\mu_1(K)h_n
  \int_{V_2^{1/3}}^\infty{{\,\d t}\over{t^2}} \cr 
&=3\mu_1(K)V_2^{2/3}h_n=3\mu_1(K)V_2^{2/3}h_0n^{-1/3}. \cr}$$
Using this estimate and Lemma 1, we get (3.8).
$\square$

\smallskip
{\bf Corollary.} \sl Let the conditions of Theorem 4 be satisfied.
Then for each $n$, 
$$\min_{h>0}\sup_x\MSE(p_n(x))
\le 3\Bigl({9\over{4\pi^2}}\Bigr)^{1/3}
\mu_1^{2/3}(K)V_2^{4/9}B^{2/3}(K)a^{2/3}n^{-2/3}.$$
\rm

\smallskip
Next we consider the so-called non-smooth case. This means that 
the underlying density function is not supposed to be differentiable 
or even continuous. Some minimum regularity conditions must be 
introduced, however (otherwise nothing substantial can be derived). 
Here this minimum condition will be the boundedness of the total
variation of the underlying density. Note that this condition
is a little less restrictive than those usually assumed 
when authors work with the non-smooth case 
(see for example van Eeden, 1985 and van Es, 1997). 

\smallskip
{\bf Theorem 5.} \sl Let the underlying density $p$ 
be a function of bounded variation, $V=\V(p)<\infty$. 
If $h_n=h_0/(\sqrt{n}\log n)$, then
$$\eqalign{
\MISE(p_{n,h})
&\le
{{\log^2n}\over{\sqrt n}}
  \Big[
  {{4\sqrt 2}\over\pi}\max\{\sqrt{\mu_1(K)},\mu_1(K)\} \cr 
&\qquad \times\max\{V^{3/2},V^2\}\,
  \max\{\sqrt{h_0},h_0\}
  +{{R(K)}\over{h_0\log n}}\Big] \cr}\eqno(3.9)$$
for all $n\ge 16$. 
\rm

\smallskip
{\bf Proof.} Let us use Lemma 2. For the second term in the 
square brackets, due to the Parseval--Plancherel identity, we have
$${1\over{nh_n}}\int |\varphi(t)|^2\,\d t
={{2\pi}\over{nh_n}}\int  K^2(x)\,\d x
={{2\pi R(K)}\over{nh_n}}.\eqno(3.10)$$
Let us estimate the first term. First we establish the following inequality:
for any $0<\alpha<1$,
$$|1-\varphi(t)|\le\mu_1(K)^\alpha 2^{1-\alpha}|t|^\alpha\eqno(3.11)$$
for all real $t$. Indeed, due to Lemma 3,
$$|1-\varphi(t)|\le\mu_1(K)|t|.\eqno(3.12)$$
For $|t|\le 2/\mu_1(K)$, the right hand side of (3.11) majorises
the right hand side of (3.12), therefore (3.11) holds for these $t$.
If $|t|>2/\mu_1(K)$, then (3.11) is evident because its right hand
side exceeds 2.

Let $\alpha$ be arbitrary inside $(0,\half)$. 
Making use of (3.11) and Lemma 4, we get
$$\eqalign{
\int |f(t)|^2|1-\varphi(h_nt)|^2\,\d t
&=2\int_0^V|f(t)|^2|1-\varphi(h_nt)|^2\,\d t \cr 
&\qquad\qquad 
  +2\int_V^\infty|f(t)|^2|1-\varphi(h_nt)|^2\,\d t \cr 
&\le 2\mu_1^{2\alpha}(K)2^{2(1-\alpha)}h_n^{2\alpha}
  \Bigl(
  \int_0^Vt^{2\alpha}\,\d t+V^2\int_V^\infty t^{2\alpha-2}\,\d t
  \Bigr) \cr
&={{2^{4-2\alpha}}\over{1-4\alpha^2}}\, 
  \mu_1(K)^{2\alpha}V^{2\alpha+1}h_n^{2\alpha}. \cr}$$
From this estimate and (3.10), using Lemma 1, we obtain
$$\eqalign{
\MISE(p_{n,h})
&\le{{2^{3-2\alpha}}\over{\pi(1-4\alpha^2)}}
 \mu_1^{2\alpha}(K)V^{2\alpha+1}h_n^{2\alpha}+{{R(K)}\over{nh_n}} \cr 
&={{2^{3-2\alpha}}\over{\pi(1-4\alpha^2)}}
  \mu_1^{2\alpha}(K)V^{2\alpha+1}h_0^{2\alpha}
  \Bigl({1\over{\sqrt n\log n}}
  \Bigr)^{2\alpha}+{{R(K)}\over{h_0}}{{\log n}\over{\sqrt n}} \cr}
  \eqno(3.13)$$
for any $\alpha\in(0,\half)$. Put
$$\alpha={{\log n}\over{2(\log n+2\log\log n)}}.$$
Then ${1\over 4}<\alpha<\half$ 
(provided that $n\ge e^e$, which translates into $n\ge16$), and hence
$$2^{3-2\alpha}<4\sqrt 2,
  \quad \mu_1(K)\le\max\{\sqrt{\mu_1(K)},\mu_1(K)\},$$
$$V^{2\alpha+1}\le\max\{V^{3/2},V^2\},
  \quad h_0^{2\alpha}\le\max\{\sqrt{h_0},h_0\}. $$
Therefore from (3.13) we obtain
$$\eqalign{
\MISE(p_{n,h})
\le&{{4\sqrt 2}\over\pi}\max\{\sqrt{\mu_1(K)},\mu_1(K)\}\,
  \max\{V^{3/2},V^2\}\, 
  \max\{\sqrt{h_0},h_0\} \cr 
&\times{1\over{1-4\alpha^2}}
  \Bigl({1\over{\sqrt n\log n}}\Bigr)^{2\alpha}
  +{{R(K)}\over{h_0}}{{\log n}\over{\sqrt n}}. \cr} \eqno(3.14)$$
Putting now 
$$\alpha_0=
{{\log n-2\log\log n}\over{2(\log n+2\log\log n)}},$$
then $\alpha>\alpha_0$ (if $n\ge e^e$), hence
$$\Bigl({1\over{\sqrt n\log n}}\Bigr)^{2\alpha}
  <\Bigl({1\over{\sqrt n\log n}}\Bigr)^{2\alpha_0}
  ={{\log n}\over{\sqrt n}}.
\eqno(3.15)$$
It remains to assess the size of $1/(1-4\alpha^2)$. We have
$$\eqalign{
{1\over{1-4\alpha^2}}
&={{(\log n+2\log\log n)^2}\over{(\log n+2\log\log n)^2-(\log n)^2}}
={{(\log n+2\log\log n)^2}\over{(2\log n+2\log\log n)2\log\log n}} \cr 
&\le{1\over4}\log n+1+{{\log\log n}\over{\log n+\log\log n}}\le\log n \cr}
  \eqno(3.16)$$ 
if $n\ge e^e$.

From (3.14), (3.15) and (3.16) we finally obtain (3.9).
$\square$

\smallskip
{\bf Corollary.} \sl Let $p$ be a unimodal density function, and 
bounded by $a$. If $h_n=h_0/(\sqrt{n}\log n)$, then
$$\eqalign{
\MISE(p_{n,h})
\le&{{\log^2n}\over{\sqrt n}}
  \cdot\Big[{{4\sqrt 2}\over\pi}\max\{\sqrt{\mu_1(K)},\mu_1(K)\} \cr 
&\times\max\{2\sqrt 2a^{3/2},a^2\}\,
  \max\{\sqrt{h_0},h_0\}+{{R(K)}\over{h_0\log n}}\Big]. \cr}$$
\rm

\section 
\centerline{\bf 4. The sinc kernel density estimator}

\start 
The sinc kernel is the function
$$K(x)={{\sin x}\over{\pi x}}$$
with the Fourier transform (defined as the principal value
of the corresponding integral)
$$\varphi(t)=\cases{
   1&for $|t|\le 1$,\cr
   0&for $|t|>1$.\cr}$$
(Sometimes the sinc kernel is defined as $K(x)=\sin(\pi x)/(\pi x)$
with the Fourier transform $\varphi(t)=I\{|t|\le\pi\}$. 
Both functions $\sin x/(\pi x)$  and $\sin(\pi x)/(\pi x)$
integrate to one in the sense of the principal value, and 
the difference is only in the scale parameter.)

From now on we focus on the kernel estimator $p_n(x)$ 
of (2.1) with $K$ being the sinc kernel. It often leads to 
better performance, and some of its properties are in fact 
easier to study than for other kernel estimators; 
see Glad, Hjort and Ushakov (1999a). Its defects -- possible 
negativeness and nonintegrability -- can easily be corrected 
by a certain modification procedure (Glad, Hjort and Ushakov, 1999b). 
It consists in setting
$$\bar p_n(x)=\max\{p_n(x)-\xi,0\}, $$
where the random $\xi$ is chosen so that the integral is 1. 
After this correction procedure, 
estimation precision of the estimator is guaranteed to improve. 

In terms of the empirical characteristic function $f_n(t)$ 
the sinc estimator can be expressed as 
$$p_n(x)={1\over{2\pi}}\int_{-1/h_n}^{1/h_n}e^{-itx}f_n(t)\,\d t.
  \eqno(4.1)$$
Suppose that the characteristic function $f$ of the 
underlying density $p$ is integrable.  
First we obtain relations for the sinc estimator, 
analogous to those of Lemmas 1 and 2 
(these cannot be applied directly since now $K$ is not integrable).

\smallskip
{\bf Lemma 5.} \sl For the sinc kernel estimator,
$$\sup_x\MSE(p_n(x))\le\Bigl\{{1\over{2\pi}}\int_{|t|\ge1/h_n}
|f(t)|\,\d t\Bigr\}^2
+{2\over{\pi nh_n}}{1\over{2\pi}}\int |f(t)|\,\d t
  \eqno(4.2)$$
and
$$\MISE(p_n)\le{1\over{2\pi}}\Bigl\{\int_{|t|\ge1/h_n}|f(t)|^2\,\d t
+{2\over{nh_n}}\Bigr\}.\eqno(4.3)$$
\rm

\smallskip
{\bf Proof.} We first prove the first inequality. We have
$$\eqalign{
\MSE(p_n&(x))
 =\E \Bigl[{1\over{2\pi}}\Bigl\{
  \int  e^{-itx}f(t)\,\d t
  -\int_{-1/h_n}^{1/h_n}e^{-itx}f_n(t)\,\d t
  \Bigr\}\Bigr]^2 \cr 
&=\E \Bigl[{1\over{2\pi}}\int_{|t|\ge1/h_n}e^{-itx}f(t)\,\d t
  +{1\over{2\pi}}\int_{-1/h_n}^{1/h_n}e^{-itx}\{f(t)-f_n(t)\}\,\d t
  \Bigr]^2 \cr 
&=\Bigl\{{1\over{2\pi}}\int_{|t|\ge1/h_n}e^{-itx}f(t)\,\d t
  \Bigr\}^2
  +\E \Bigl[
  {1\over{2\pi}}\int_{-1/h_n}^{1/h_n}e^{-itx}\{f(t)-f_n(t)\}\,\d t
  \Bigr]^2. \cr}$$
Let us estimate the second term on the right hand side. Denote it
by $T_2$. Taking into account that 
$$\E f_n(u)f_n(v)=(1-1/n)f(u)f(v)+(1/n)f(u+v), $$
we obtain
$$\eqalign{
T_2&={1\over n}{1\over{(2\pi)^2}}
  \int_{-1/h_n}^{1/h_n}\int_{-1/h_n}^{1/h_n}
  e^{-i(u+v)x}\{f(u+v)-f(u)f(v)\}\,\d u\,\d v \cr 
&={1\over{2\pi n}}\int_{-1/h_n}^{1/h_n}\Bigl\{
  {1\over{2\pi}}\int_{-1/h_n}^{1/h_n}e^{-i(u+v)x}f(u+v)\,\d u
  \Bigr\}\,\d v \cr 
&\qquad\qquad\qquad 
  -{1\over n}\Bigl\{{1\over{2\pi}}\int_{-1/h_n}^{1/h_n}e^{-itx}f(t)\,\d t
  \Bigr\}^2 \cr 
&={1\over{2\pi n}}\int_{-1/h_n}^{1/h_n}\Bigl\{
  {1\over{2\pi}}\int_{-1/h_n+v}^{1/h_n+v}e^{-itx}f(t)\,\d t
  \Bigr\}\,\d v
  -{1\over n}\Bigl\{{1\over{2\pi}}\int_{-1/h_n}^{1/h_n}e^{-itx}f(t)\,\d t
  \Bigr\}^2. \cr}$$
Therefore
$$\eqalign{
T_2&\le{1\over{2\pi n}}\int_{-1/h_n}^{1/h_n}\,\d v\,
    {1\over{2\pi}}\int |f(t)|\,\d t
   +{1\over n}{1\over{2\pi}}\int |f(t)|\,\d t\, 
    {1\over{2\pi}}\int_{-1/h_n}^{1/h_n}\d s \cr 
&={2\over{\pi nh_n}}{1\over{2\pi}}\int |f(t)|\,\d t. \cr}$$
Thus we obtain (4.2).

Next we prove (4.3). Observe that we may use relation (2.5) with
$$\hat f_n(x)=\cases{
  f_n(x)&if $|t|\le1/h_n$,\cr
  0& otherwise. \cr}$$
Therefore
$$\MISE(p_n)={1\over{2\pi}}\Bigl\{
  \int_{-1/h_n}^{1/h_n}\E|f_n(t)-f(t)|^2\,\d t+
  \int_{|t|\ge1/h_n}|f(t)|^2\,\d t\Bigr\},$$
and it suffices to show that
$$\int_{-1/h_n}^{1/h_n}\E |f_n(t)-f(t)|^2\,\d t\le{2\over{nh_n}}.$$
Taking now into account that
$\E |f_n(t)|^2=(1-1/n)|f(t)|^2+1/n$, we obtain
$$\int_{-1/h_n}^{1/h_n}\E |f_n(t)-f(t)|^2\,\d t={1\over n}
  \int_{-1/h_n}^{1/h_n}(1-|f(t)^2|)\,\d t\le{1\over n}
  \int_{-1/h_n}^{1/h_n}\,\d t={2\over{nh_n}}.$$
This proves the claim. $\square$

\smallskip
Now we derive some estimates for $\MISE$ and $\MSE$ of the sinc
estimator in terms of the degree of smoothness 
of the underlying density. First we consider 
the non-smooth case, when a density to be estimated 
is not supposed to be differentiable or even continuous.

\smallskip
{\bf Theorem 6.} \sl Let $p$ have bounded variation, 
$\V(p)=V<\infty$, and let $p_n$ be the sinc estimator. 
If $h_n=h_0/\sqrt n$, then
$$\MISE(p_n)\le
{1\over{\pi\sqrt n}}\Bigl(V^2h_0+{1\over{h_0}}\Bigr).$$
\rm

\smallskip
{\bf Proof.} Making use of relation (4.3) of Lemma 5 
and Lemma 4, we obtain
$$\eqalign{
\MISE(p_n)
&\le {1\over{2\pi}}\Bigl\{
  \int_{|t|\ge 1/h_n}|f(t)|^2\,\d t+{2\over{nh_n}}\Bigr\} \cr 
&\le{1\over{2\pi}}\Bigl(2V^2\int_{1/h_n}^\infty{{\,\d t}\over{t^2}}+
  {2\over{nh_n}}\Bigr)
  ={1\over{\pi\sqrt n}}\Bigl(V^2h_0+{1\over{h_0}}\Bigr), \cr}$$
as required. $\square$
\eject 

\smallskip
{\bf Corollary 1.} \sl Let the conditions of Theorem 6 be satisfied.
Then for each $n$,
$$\min_{h>0}\MISE(p_n)\le{{2V}\over{\pi\sqrt n}}.$$
\rm

\smallskip
{\bf Corollary 2.} \sl Let $p$ be a unimodal density function,
and let $p_n$ be the sinc estimator. If $p$ is bounded by $a$, 
and $h_n=h_0/\sqrt n$, then
$$\MISE(p_n)\le
{1\over{\pi\sqrt n}}\Bigl(4a^2h_0+{1\over{h_0}}\Bigr),$$
and
$$\min_{h>0}\MISE(p_n)\le{{4a}\over{\pi\sqrt n}}.$$
\rm

\smallskip
Now consider the case when the density to be estimated is $m$ times
differentiable, $m\ge 1$. It will be shown that in this case the
upper bound for $\MISE$ of the sinc estimator has order $n^{-2m/(2m+1)}$
that in principal cannot be achieved (for $m>2$) for 
kernel estimators with kernels being density functions.

\smallskip
{\bf Theorem 7.} \sl Let $p$ be $m$ times differentiable
with $p^{(m)}$ a function of bounded variation, 
$\V(p^{(m)})=V_m<\infty$. If $p_n$ is the sinc estimator, and
$h_n=h_0n^{-1/(2m+1)}$, then
$$\MISE(p_n)\le{1\over{2\pi}}
  \Bigl\{{4(m+1)\over{2m+1}}V_m^{(2m+1)/(m+1)}h_0^{2m}+{2\over{h_0}}
  \Bigr\}n^{-2m/(2m+1)}. \eqno(4.4)$$
\rm

\smallskip
{\bf Proof.} We have
$$\eqalign{
\int_{|t|\ge 1/h_n}|f(t)|^2\,\d t
&=h_n^{2m}\int_{|t|\ge 1/h_n}
  \Bigl({1\over{h_n}}\Bigr)^{2m}|f(t)|^2\,\d t \cr 
&\le h_n^{2m}\int_{|t|\ge 1/h_n}|t|^{2m}|f(t)|^2\,\d t
  \le h_n^{2m}\int |t|^{2m}|f(t)|^2\,\d t. \cr}\eqno(4.5)$$
Let us estimate the integral on the right hand side, making use of
Lemma 4. We have
$$|f(t)|=\cases{
1&for $|t|\le V_m^{1/(m+1)}$,\cr
V_m/|t|^{m+1}&for $|t|>V_m^{1/(m+1)}$,\cr}$$
therefore
$$\eqalign{
\int |t|^{2m}|f(t)|^2\,\d t
&\le 2\int_0^{V_m^{1/(m+1)}}t^{2m}\,\d t+
  2V_m^2\int_{V_m^{1/(m+1)}}^\infty{{\,\d t}\over{t^2}} \cr 
&={{4(m+1)}\over{2m+1}}V_m^{(2m+1)/(m+1)}. \cr}\eqno(4.6)$$
Thus, from inequality (4.3) of Lemma 5, and relations 
(4.5) and (4.6), we obtain (4.4).
$\square$

\smallskip
{\bf Corollary.} \sl Let the conditions of Theorem 7 be satisfied.
Then for each $n$,
$$\min_{h>0}\MISE(p_n)
  \le{1\over{2\pi}}\{4(m+1)\}^{1/(2m+1)}
  \Bigl({{2m+1}\over m}\Bigr)^{2m/(2m+1)}V_m^{1/(m+1)}n^{-2m/(2m+1)}.$$
\rm

\smallskip
{\bf Theorem 8.} \sl Let $p$ be $m$ times differentiable, 
with $p^{(m)}$ a function of bounded variation, 
$\V(p^{(m)})=V_m<\infty$. If $p_n$ is the sinc estimator, 
and $h_n=h_0n^{-1/(2m-1)}$, then
$$\eqalign{
\sup_x\MSE(p_n(x))
\le&{1\over{\pi^2}}
  \Big\{{{(m+1)^2}\over{m^2}}V_m^{2m/(m+1)}h_0^{2(m-1)} \cr 
&+2\Big(V_m^{1/(m+1)}+{1\over m}V_m^{m/(m+1)}\Big){1\over{h_0}}
  \Big\}n^{-2(m-1)/(2m-1)}. \cr}$$
\rm

\smallskip
The proof of this theorem is analogous to that of Theorem 7, 
one just needs to use relation (4.2) of Lemma 5 instead of relation 
(4.3) and take into account that due to Lemma 4,
$$\eqalign{A(p)
 ={1\over{2\pi}}\int |f(t)|\,\d t
  &\le{1\over{2\pi}}\int_{-V_m^{1/(m+1)}}^{V_m^{1/(m+1)}}\,\d t
  +{1\over{2\pi}}\int_{|t|>V_m^{1/(m+1)}}{V_m\,\d t\over{|t|^{m+1}}} \cr 
&\le{1\over\pi}
  \Bigl\{V_m^{1/(m+1)}+{1\over m}V_m^{m/(m+1)}\Bigr\}. \cr}$$

\smallskip
{\bf Corollary.} \sl Let the conditions of Theorem 8 be satisfied.
Then for each $n$,
$$\eqalign{
\min_{h>0}\sup_x\MSE(p_n(x))
\le &{{2m-1}\over{\pi^2}}\Bigl({{m+1}\over m}\Bigr)^{2/(2m-1)} \cr 
&\times
  \Big\{{{m+V_m^{(m-1)/(m+1)}}\over{m(m-1)}}\Big\}^{{2m-2}\over{2m-1}}
  V_m^{2/(m+1)}n^{-2(m-1)/(2m-1)}. \cr}$$
\rm

\smallskip
Now we proceed to the `supersmooth' case which we define in terms
of characteristic functions (although this class of distribution can
be defined in terms of density functions as well, a description in
terms of characteristic functions is simpler, more natural
and more convenient for our purposes). A distribution $F$ with
characteristic function $f(t)$ is said to be supersmooth if
for some $\alpha>0$ and $\gamma>0$,
$$B(p;\alpha,\gamma)=\int e^{\gamma|t|^\alpha}|f(t)|\,\d t<\infty.$$
Thus a normal density is supersmooth with $\alpha=2$
while a Cauchy is supersmooth with $\alpha=1$, for example. 

\smallskip
{\bf Theorem 9.} \sl Let the characteristic function $f$ of
$p$ have a finite $B(p;\alpha,\gamma)$ value, 
for some positive $\alpha$ and $\gamma$. 
If $p_n$ is the sinc estimator, and
$$h_n=\Bigl\{{1\over\gamma}\log(h_0n)\Bigr\}^{-1/\alpha},
  \quad {\sl with\ } h_0\ge{1\over n},$$
then
$$\MISE(p_n)\le{1\over{2\pi n}}
\Bigl\{
  {2\over{\gamma^{1/\alpha}}}(\log h_0+\log n)^{1/\alpha}
  +{{B(p;\alpha,\gamma)}\over{h_0}}
  \Bigr\}.
\eqno(4.7)$$
\eject 
\rm

\smallskip
{\bf Proof.} We have 
$$\eqalign{
\int_{|t|\ge 1/h_n}|f(t)|^2\,\d t
&\le\int_{|t|\ge 1/h_n}|f(t)|\,\d t
  =e^{-\gamma/h_n^\alpha}\int_{|t|\ge 1/h_n}
   e^{\gamma/h_n^\alpha}|f(t)|\,\d t \cr 
&\le e^{-\gamma/h_n^\alpha}
  \int  e^{\gamma|t|^\alpha}|f(t)|\,\d t
  ={{B(p;\alpha,\gamma)}\over{nh_0}}. \cr}$$
Using this estimate and inequality (4.3) of Lemma 5, we obtain (4.7).
$\square$

\smallskip
{\bf Theorem 10.} \sl Let the conditions of Theorem 9 be satisfied,
and let again $A(p)=(2\pi)^{-1}\int|f(t)|\,\d t$. Then
$$\sup_x\MSE(p_n(x))\le{1\over{n}}
\Bigl\{
{2A(p)\over{\pi\gamma^{1/\alpha}}}(\log h_0+\log n)^{1/\alpha}
  +{{B^2(p;\alpha,\gamma)}\over{4\pi^2nh_0}}\Bigr\}.$$
\rm

\smallskip
The proof of the theorem is similar to that of Theorem 9 
(inequality (4.2) is used instead of (4.3)).

Theorems 9 and 10 can be improved for one
subclass of supersmooth densities. The result is given by the
next theorem and is quite curious. Note that this theorem corresponds
to a special case of a result by Ibragimov and Khas'minskii (1982),
and we give it here for the completeness.

\smallskip
{\bf Theorem 11.} \sl Let the characteristic function $f$ of
$p$ satisfy the condition: there exists $\tau>0$ such that
$f(t)=0$ for $|t|>\tau$. If $p_n$ is the sinc estimator, and
$h_n\le1/\tau$, then
$$\sup_x\MSE(p_n(x))\le{{2A(p)}\over{\pi nh_n}}
\le{{2\tau}\over{\pi^2nh_n}},$$
and
$$\MISE(p_n)\le{1\over{\pi nh_n}}.$$
In particular, if $h_n={\rm const.}=1/\tau$, then
$$\sup_x\MSE(p_n(x)) \le{{2\tau^2}\over{\pi^2n}}
  \quad {\rm and} \quad 
  \MISE(p_n)\le{\tau \over{\pi n}}.$$
\rm

\smallskip
A proof of the theorem can be immediately obtained from inequalities 
(4.2) and (4.3) of Lemma 5:
integrals on the right hand sides of these vanish when $h_n\le 1/\tau$.

Theorem 11 implies in particular that if the characteristic 
function of the underlying distribution vanishes for large
values of the argument, and one uses the sinc estimator for
approximation, then $p_n$ converges to $p$ as $n\to\infty$
even when $h_n$ does not converge to zero.
\eject 

\section 
\centerline{\bf 5. Discussion and applications}  

\start
This article has provided upper bounds for
both the traditional MISE and also the less worked with 
max-MSE performance measures of kernel density estimators. 
A list of such upper bounds has been provided, under 
various sets of assumptions, for both the traditional
kernels as well as for the sinc kernel, which has 
particularly attractive features. 
Our finite-sample results have been reached entirely outside 
the customary framework of asymptotics, Taylor expansions
and small bandwidths, through the extensive use of 
characteristic and empirical characteristic functions. 
Below we give some concluding remarks, pointing to
ways in which the results can be applied in statistics. 

\smallskip
{\sl 5.1. Rule-of-thumb bandwidths for MISE and max-MSE.}  
Consider kernel estimators with a traditional 
kernel $K$, a symmetric density. 
The traditional large-sample approximations 
lead to an asymptotically optimal bandwidth of size 
$$h_n=\{R(K)/\mu_2^2(K)\}^{1/5}R(p'')^{-1/5}\,n^{-1/5}, $$
and with consequent minimum approximate MISE of size 
$$\min_{h>0}{\rm AMISE}(p_n)=(5/4)\{\mu_2^2(K)R^4(K)\}^{1/5}
  R(p'')^{1/5}\,n^{-4/5}, $$
see for example Wand and Jones (1995). 
When $K$ is standard normal, and $p$ is a normal density with 
standard deviation $\sigma$, this leads to the popular 
`normal rule-of-thumb' bandwidth $h_n=1.0592\,\sigma\,n^{-1/5}$. 
Note that the structure of these classic results is 
very similar to that seen in Theorem 1 and its corollary;
in particular, the well-known large-sample result about 
the $n^{-4/5}$ precision rate is here reached entirely 
without asymptotics machinery or approximations. 

It is interesting to compare the above with what one finds 
using the upper bounds. For a normal density, 
$V_2=\int|p'''|\,\d x$ is found to be the scale factor $\sigma^{-3}$
times $2(2\pi)^{-1/2}\{1+4\exp(-3/2)\}=1.5100$,
and this leads via Theorem 1 to the rule $h_n=0.8204\,\sigma\,n^{-1/5}$. 
This has been calculated using upper bound results 
derived under minimal assumptions, and 
which hence do not pretend to be very accurate for smooth densities
like the normal. It is comforting to see that only a moderate
amount is lost in precision, in this very smooth case, 
since the ratio of the minimised upper bound 
to the minimised asymptotic MISE is found to be 1.2911. 

The max-MSE criterion is a natural venue, seemingly not travelled before.
It is difficult to reach applicable results for this criterion
based on the traditional approximations. However, 
Theorem 3 and its corollary provide ways of bounding the max-MSE 
when there is information on $V_3=\int|p^{(4)}|\,\d x$. 
If $p$ is normal, then $V_3$ can be shown to be 
the scale factor $\sigma^{-4}$ times 
$4\{(3b-b^3)\phi(b) - (3c-c^3)\phi(c)\}=2.8006$, where 
$b=(3-\sqrt{6})^{1/2}$ and $c=(3+\sqrt{6})^{1/2}$,
and $\phi$ is the standard normal density. 
The normal rule-of-thumb, when the normal kernel is used,
becomes $h_n=1.1883\,\sigma\,n^{-1/5}$, 
which again is not far from the traditional rule-of-thumb. 

We also point out that some of these results may 
be sharpened under further assumed constraints on
the underlying density. The quite crude bound (3.2)
has for example been used for $|f(t)|$, which 
could be bounded more effectively under such additional
restrictions. This is not pursued here, however. 

\smallskip 
{\sl 5.2. Cross-validation and normal rule-of-thumb in new light.} 
Results reached in this article, about MISE and upper bounds expressed 
in terms of characteristic functions, point to new ways in 
which bandwidths can be selected from data. 
Expressions (2.8) and its upper bound given in Lemma 2 
depend on $q(t)=|f(t)|^2$, but not on other aspects of 
the underlying density $p$. A suitable estimate $\hat q(t)$ 
may now be inserted in these expressions, after which 
one may minimise over the smoothing parameter $h$. 
For the normal kernel, this could for example mean minimising 
$$Q_n(h)={1\over 2\pi}\int \hat q(t)\{1-\exp(-\half h^2t^2)\}^2\,\d t
  +{1\over 2\pi}{1\over n}\int \exp(-h^2t^2)\,\d t $$
over $h$, after having selected an estimator $\hat q(t)$. 
Interestingly it turns out that this scheme reproduces 
the well-known `unbiased cross-validation' rule, 
see for example Wand and Jones (1995), 
when one employs the natural nonparametric unbiased estimator  
$$\hat q(t)={1\over n(n-1)}\sum_{j\not=k}\exp(it(X_j-X_k))
  ={2\over n(n-1)}\sum_{j<k}2\cos(t(X_j-X_k)). $$

Other methods emerge by using alternative estimators for $q(t)$.
If one above uses the simplest parametric estimate, 
namely $\exp(-\hat\sigma^2t^2)$ (or a debiased version thereof) 
corresponding to a normality approximation, 
then minimising $Q_n(h)$ is a better finite-sample version
of the classic rule-of-thumb $1.0592\,\hat\sigma\,n^{-1/5}$.
Semiparametric versions of this argument would be worth
studying; see Hjort (1999) for a similar enterprise. 

\smallskip 
{\sl 5.3. Minimax precision control.} 
Another type of application would be following.
Assume that an upper bound for $V_2=\int |p'''|\,\d x$
is established, say $V_2\le v_2$. This is a statement
of the maximum envisaged wigglyness of the density;
a small $V_2$ would mean a density which can be approximated
with a quadratic function. The bound $v_2$ could 
be set after inspection of data or from prior grounds. -- In such
a situation the corollary of Theorem 1 can be given to
find a sample size $n_0$ guaranteed to secure a MISE-accuracy
below a given threshold, for the big class of all densities 
with $V_2\le v_2$. 

Variations of this scenario can easily be given, 
for example selecting a sample size necessary to secure 
max-MSE below a certain precision threshold for all 
relevant densities, as constrained with bounds 
on $V_3$ and the maximum value $a$. 

\smallskip
{\sl 5.4. Large-sample superiority of the sinc method.} 
Approximation of densities via the sinc kernel can often
be more accurate than with traditional kernels. 
This has been pointed out early on by Davis (1977),
but the method does not seem popular in practice. 
That it takes negative values and is not (Lebesgue)-integrable
is not a real concern, since it can be repaired for these defects
in an automatic fashion which also guarantees 
precision improvement; see Glad, Hjort and Ushakov (1999b). 
And results from Section 4 above, with further analysis
provided in Glad, Hjort and Ushakov (1999a), give
clear indications of the strong performance of the sinc method. 

Theorem 9 shows that the rate of MISE towards zero 
is often much better than the $n^{-4/5}$ available 
with traditional kernels. For any mixture of normals,
for example, the rate is $O((\log n)^{1/2}/n)$,
while if some Cauchy type tail behaviour is mixed in
it becomes $O((\log n)/n)$. The same remarks apply 
to the max-MSE performance criterion. 

That the sinc method also can perform better than 
the traditional ones for non-smooth cases is made 
clear by Theorem 6, where a $O(n^{-1/2})$ rate is exhibited
for MISE under a minimal assumption on $p$, 
compared to the $O(n^{-1/2}\log^2n)$ rate for the ordinary methods. 

\bigskip\noindent
{\bf Acknowledgments.} 
The second author is grateful for the opportunity to
visit the Department of Mathematics at the University of Oslo
during the summer of 1998, and the first author equally grateful 
for being able to visit 
{\cyr Mehaniko-matematicheski{\j3} fakulp1tet 
Moskovskogo gosudarstvennogo universiteta} in 1999.
These visits have been made possible with 
partial funding from the Norwegian Research Foundation
and the Department of Mathematics at the University of Oslo. 
The authors also appreciate fruitful discussions with Ingrid Glad. 

\def\ref#1{{\noindent\hangafter=1\hangindent=20pt
  #1\smallskip}}
\parindent0pt
\baselineskip11pt
\parskip3pt

\bigskip
\centerline{\bf References}

\medskip
\ref{%
Davis, K.B. (1977).
Mean integrated square error properties of density estimates. 
{\sl Annals of Statistics} {\bf 5}, 530--535.}

\ref{%
Devroye, L.~and Gy\"orfi, L. (1985). 
{\sl Nonparametric Density Estimation: The $L_1$ View.} 
Wiley, New York.}

\ref{%
van Eeden, C. (1985). 
Mean integrated squared error of kernel estimators 
when the density and its derivatives are not necessarily continuous. 
{\sl Annals of the Institute of Statistical Mathematics} {\bf 37}, 
Part A, 461--472.}

\ref{%
van Es, A.J. (1997). 
A note on the integrated squared error of a 
kernel density estimator in non-smooth cases. 
{\sl Statistics and Probability Letters} {\bf 35}, 241--250.}

\ref{%
Feller, W. (1971). 
{\sl An Introduction to Probability Theory
and its Applications}, Vol.~2 (2nd ed.). Wiley, New York.}
\eject 

\ref{%
Glad, I.K., Hjort, N.L.~and Ushakov, N.G. (1999a). 
Density estimation using the sinc kernel.
{\sl Statistical Research Report},  
Department of Mathematics, University of Oslo.}

\ref{%
Glad, I.K., Hjort, N.L.~and Ushakov, N.G. (1999b). 
Correction of density estimators which are not densities. 
{\sl Statistical Research Report},  
Department of Mathematics, University of Oslo.}

\ref{%
Hjort, N.L. (1999).
Towards semiparametric bandwidth selectors for kernel density estimators.
{\sl Statistical Research Report},  
Department of Mathematics, University of Oslo.}

\ref{%
Ibragimov, I.A.~and Khas'minskii, R.Z. (1982). 
Estimation of distribution density belonging to a class of entire functions.
{\sl Theory of Probability and their Applications} {\bf 27}, 
No. 3, 551--562.}

\ref{%
Titchmarsh, E. (1937). 
{\sl Introduction to the Theory of Fourier Integrals.} 
Clarendon Press, Oxford.}

\ref{%
Ushakov, N.G.~and Ushakov, V.G. (1999). 
Some inequalities for characteristic functions 
of densities with bounded variation. 
{\sl Statistical Research Report},  
Department of Mathematics, University of Oslo.}

\ref{%
Wand, M.P.~and Jones, M.C. (1995).
{\sl Kernel Smoothing.} Chapman and Hall, London.} 

\bigskip
\bigskip
\baselineskip12pt 
\settabs 2\columns 
\+\qquad Nils Lid Hjort            &Nikolai G.~Ushakov \cr 
\+\qquad Department of Mathematics &Institute of Microelectronics Technology\cr
\+\qquad University of Oslo        &{\cyr Rossi{\j3}skaya Akademiya Nauk} \cr 
\+\qquad P.B.~1053 Blindern        &{\cyr 142 432 Chernogolovka} \cr 
\+\qquad Oslo, Norway              &{\cyr Moskovskaya oblastp1, Rossiya} \cr

\baselineskip10pt 
\bigskip
\line{\hskip7.12truecm\smallrm 1999 \hfill} 
\line{\hskip7.12truecm\smallrm 2000 \hfill} 
\line{\hskip7.12truecm\smallrm 3111 \hfill} 
\vskip-15pt 
$$\vdots$$

\bye